\documentclass[12pt]{article}
\usepackage{amssymb}
\def\ba#1{\begin{array}{#1}}
\def\ea{\end{array}}
\def\beq#1{\begin{equation}\label{#1}}
\def\eeq{\end{equation}}

\newcommand{\RR}{\mathbb{R}}

\newcommand{\CC}{\mathbb{C}}

\newcommand{\HH}{\mathbb{H}} 
\newcommand{\cM}{\mathcal{M}} 
\newcommand{\cP}{\mathcal{P}} 
\newcommand{\mS}{\hat{\St}} 
\newcommand{\mL}{\hat{\Lt}} 
\newcommand{\mP}{\hat{\cP}} 
\newcommand{\Cp}{\CC_+}
\newcommand{\St}{\mathcal{S}} 
\newcommand{\Ht}{\mathcal{H}} 

\newcommand{\f}{\frac}
\newcommand{\unl}{\underline}

\newtheorem{thm}{Theorem}
\newtheorem{lem}{Lemma}

\newcommand{\ph}{\varphi}

\newcommand{\Lp}{L^{p}}

\newcommand{\Lt}{\mathcal{L}} 
\newcommand{\Lta}{\Lt_\theta}

\newcommand{\essup}{\mathop{\mathrm{ess\, sup}}}
\def\dst{\displaystyle}
\def\unl{\underline}
\newcommand{\pf}{\noindent{\it Proof}. }
\def\eop{\hfill$\Box$}
\def\rem#1{}

\begin{document}
\begin{center}
{\Large
$\Lp$ estimates for angular maximal functions
associated with
\\[0.5ex]
 Stieltjes and Laplace transforms}
\\[3ex]
{\bf Sergey Sadov}
\\[3ex]
{\small
Memorial University of Newfoundland, Canada}
\\[1ex]
{\small sergey@mun.ca}
\end{center}

\bigskip\noindent
{\small
{\bf Abstract}.
  Maximal angular operator sends a function defined in a sector of the complex plane with vertex at 0
  to the function of modulus obtained by
  maximizing over argument. 
  Compositions of the so defined maximal angular operator
  (in suitable sectors) with the Poisson, Stieltjes and   
  Laplace transforms are shown to be bounded (nonlinear)         
  operators from $L^p$ to $L^q$ for the same values of $p$    
  and $q$ as their standard counterparts.
}

\bigskip\noindent
{\bf 1.}\@
Suppose $g(z)$ is a complex-valued function defined in the sector
$$
\CC_{(\theta_1,\theta_2)}=\{z\in\CC\,|\,\theta_1<\arg z<\theta_2\}.
$$

The {\it angular maximal function}\ for $g$ with respect
to $\CC_{(\theta_1,\theta_2)}$ is a function $\RR_+\to [0,\infty]$ defined as follows:
\beq{fmsec}
 \cM_g^{(\theta_1,\theta_2)}(\rho)=\essup\limits_{\theta\in(\theta_1,\theta_2)} |g(\rho e^{i\theta})|.
\eeq
In this work, the function $g(z)$ will be harmonic
in the sector of interest and $\essup$
in (\ref{fmsec}) can be replaced by $\sup$.
For convenience we introduce shorter notation for the angular maximal functions corresponding to
the three sectors that will be used: (i) the plane cut along the positive real axis $\CC^*=\CC_{(0,2\pi)}$,
(ii) the upper half-plane $\HH=\CC_{(0,\pi)}$,
and (iii) the right half-plane $\Cp=\CC_{(-\pi/2,\pi/2)}$:
$$
\begin{array}{ccl}
 {\rm (i)}&\quad &\dst
\cM^*_g(\rho)=\sup_{0<\theta<2\pi} |g(\rho e^{i\theta})|,
\\[2ex]
 {\rm (ii)}&\quad &\dst
 \cM^u_g(\rho)=\sup_{0<\theta<\pi} |g(\rho e^{i\theta})|,
\\[2ex]
 {\rm (iii)}&\quad &\dst
 \cM^+_g(\rho)=\sup_{|\theta|<\pi/2} |g(\rho e^{i\theta})|.
\end{array}
$$

\bigskip\noindent
{\bf 2.}\@
The objects of our study are the angular maximal functions associated with
the Poisson integral, the Stieltjes transform, and the Laplace transform.

\smallskip
(i) If $f$ is defined on $\RR$ and $\int_{\RR} f(t)(t^2+1)^{-1}\,dt<\infty$, then the {\em Poisson integral}\
with density $f(t)$ is a harmonic function in $\HH$, 
\beq{Pt}
 \cP f(x+iy)=\f{1}{\pi}\int_{-\infty}^\infty \f{y}{(t-x)^2+y^2}\,f(t)\,dt.
\eeq

\smallskip
(ii) If $f$ is defined on $\RR_+$ and 
$\int_{\RR_+}f(t)(t+1)^{-1}\,dt<\infty\;$, then the {\em Stieltjes transform}\ of $f$ is an analytic function in $\CC^*$,
\beq{St}
 \St f(z)=\int_{0}^\infty \f{f(t)}{t-z}\,dt.
\eeq

\smallskip
(iii) The {\em Laplace transform}\ of a function $f(t)$ defined on $\RR_+$ (say, $f\in \Lp(\RR_+)$ for some $p\in[1,\infty]$)
is an analytic function in $\Cp$,
\beq{Lt}
 \Lt f(z)=\int_0^\infty f(t)e^{-zt}\,dt.
\eeq

The {\it angular maximal Poisson transform}\ $\mP$
sends a function $f(t)$ defined on $\RR$ to a function of $\rho\in\RR_+$,
\beq{mpt}
 \mP f(\rho)=\cM^u_{\cP f}(\rho)=\f{1}{\pi}
 \sup_{0<\theta<\pi} \left|\int_{-\infty}^\infty\f{\rho\sin\theta}{(t-\rho\cos\theta)^2+(\rho\sin\theta)^2}\,f(t)\,dt\right|.
\eeq

The {\it angular maximal Stieltjes transform}\ $\mS$ sends a function
$f(t)$ defined on $\RR_+$ to a function of $\rho\in\RR_+$,
$$
 \mS f(\rho)=\cM^*_{\St f}(\rho)=\sup_{0<\theta<2\pi} \left|\int_0^\infty\f{f(t)}{t-\rho e^{i\theta}}\,dt\right|.
$$

The {\it angular maximal Laplace transform}\ $\mL$ sends a function
$f(t)$ defined on $\RR_+$ to a function of $\rho\in\RR_+$,
$$
 \mL f(\rho)=\cM^+_{\Lt f}(\rho)=\sup_{|\theta|<\pi/2} \left|\int_0^\infty e^{-t\rho e^{i\theta}}\,f(t)\,dt\right|.
$$

\bigskip\noindent
{\bf 3.}\@ Our goal is to obtain $L^p$ estimates for the maximal transformations introduced above.
Here the results are formulated; proofs follow in the subsequent sections.

The key technical result is an estimate of weak type $(1,1)$  for the maximal angular Poisson transform.
We use the notation
$$
 E_f(\lambda)=\{x:\, |f(x)|>\lambda\},
 \qquad
 \mu_f(\lambda)=|E_f(\lambda)|.
$$
(Here $|\cdot|$ denotes the Lebesgue measure of a subset of $\RR$.)

\begin{thm}\label{mwPthm}
The map $f\mapsto\mP f$ is of weak type (1,1). Specifically,
if $f\in L^1(\RR)$, then for any $\lambda>0$
\beq{wL1P}
\mu_{\mP f}(\lambda)\leq K_1\f{\|f\|_1}{\lambda}
\eeq
with constant $K_1$ independent of $f$.
\end{thm}

An easy corollary, by means of the Marcinkiewicz interpolation theorem, is 

\begin{thm}\label{mPthm}
Let $f\in L^p(\RR)$, $1<p\leq\infty$. Then
\beq{LpP}
\|\mP f\|_p\leq K_2\|f\|_p,
\eeq
where $K_2$ depends only on $p$ but not on $f$.
\end{thm}

Throwing in the $L^p$ boundedness of the Hilbert transform, we obtain the following maximal theorem for the
Stieltjes transform. 

\begin{thm}\label{mSthm}
Let $f\in L^p(\RR_+)$, $1<p<\infty$. Then
\beq{LpS}
\|\mS f\|_p\leq K_3\|f\|_p,
\eeq
where $K_3$ depends only on $p$ but not on $f$.
\end{thm}

Next, a ``maximal angular theorem'' for the Laplace transform 
(mentioned in \cite[\S~4.3]{MS} as a hypothetical result)
will be derived by means of the
Cauchy formula and using an estimate for the Laplace transform along rays in the right half-plane.
We always assume that the exponents $p$ and $p'$ are conjugate, i.e.  $p^{-1}+p'^{-1}=1$.

\begin{thm}\label{mLthm}
Let $f\in L^p(\RR_+)$, $1\leq p\leq2$. Then
\beq{LpL}
\|\mL f\|_{p'}\leq K_4\|f\|_p,
\eeq
where $K_4$ depends only on $p$ but not on $f$.
\end{thm}

\noindent{\bf Remark}. These theorems can be restated in terms of uniform norm bounds
for families of linear integral operators parametrized by measurable functions $\theta(r)$ corresponding to curves $z=re^{i\theta(r)}$, 
along which the transformed functions are observed.

\smallskip
For instance, a generalization of Theorem~\ref{mSthm} 
can be obtained as its corollary. 
For a measurable function $\,\theta(x):\RR_+\to(0,\pi)$
denote 
$\St_\theta f(x)=\int_{\RR} (xe^{i\theta(x)}-y)^{-1}\,f(y)\, dy$.
By Theorem~\ref{mSthm}, $\|\St_\theta\|_{L^p(\RR)\to L^p(\RR_+)}\leq 2K_3(p)$. (The norm bound is independent of $\theta$.)\@ 
Computing the kernel of the composition 
$\St_{\theta_1}\St^*_{\theta_2}$,
we get the following result.

\begin{thm}\label{mHilb}
For any two measurable functions $\theta_1(r)$ and $\theta_2(r)$ from $\RR_+$ to
$(0,\pi)$ the integral operator 
$$
\St_{\theta_1(\cdot),\theta_2(\cdot)} f(x)=\int_{0}^\infty 
\frac{f(y)}{xe^{i\theta_1(x)}-ye^{-i\theta_2(y)}}
\, dy
$$ 
is bounded in $L^p(\RR_+)$, $1<p<\infty$. 
A norm bound $K_5(p)$ common to all pairs of functions
$\{\theta_1(\cdot),\theta_2(\cdot)\}$ exists.
\end{thm}

Consequently, the operators with kernels $(x-ye^{i(\theta_1(x)+\theta_2(y))})^{-1}$
are uniformly bounded in  $L^p(\RR_+)$, $1<p<\infty$, for all pairs of functions 
$\theta_1(x),\theta_2(y)$ from $\RR_+$ to $(0,\pi)$. This is 
the aforementioned generalization-corollary of Theorem~\ref{mSthm}. 

\smallskip
On the other hand, for any $p>1$ the norm in $L^p(\RR)$
of the operator with kernel
$$
 \mathcal{K}_{\epsilon}(x,y)=\left\{\begin{array}{l}
  (x-ye^{i\epsilon})^{-1} \;\;\;\; {\rm if}\;\; y<x\\
  (x+ye^{-i\epsilon})^{-1}\,\;\; {\rm if}\;\; y>x
\end{array}
\right.
$$
tends to $\infty$ as $\epsilon\to 0^+$. Thus, the family of operators with kernels  $(x-ye^{i\theta(x,y)})^{-1}$, where $\theta:\,\RR^2\to (0,2\pi)$ is an arbitrary
measurable function, is not iniformly bounded. 

\bigskip\noindent
{\bf 4.}\@
We begin proofs from an easy end, first demonstrating the implications \\[\smallskipamount]
\centerline{
Theorem~\ref{mwPthm} $\Longrightarrow$
Theorem~\ref{mPthm}$\Longrightarrow$ Theorem~\ref{mSthm}. 
}

\smallskip
Classical facts of harmonic analysis used in this paper --- Marcinkiewicz's theorem,
$L^p$-boundedness of the Hilbert transform, properties of the Poisson integral, the Hardy-Littlewood
maximal theorem --- can be found, for example, in \cite{Grafakos}, \cite{Stein}.

\smallskip
If $f\in L^{\infty}(\RR)$, then for any $x\in\RR$ and $y>0$ the estimate
$$
 |\cP f(x+iy)|\leq \|f\|_\infty
$$
holds due to the approximate identity properties of the Poisson kernel.
Hence $\|\mP f\|_\infty\leq \|f\|_\infty$. This estimate together with (\ref{wL1P})
immediately leads to (\ref{LpP}) due to Marcinkiewicz's theorem.

\smallskip
Now, if $z\in\HH$ and $f$ is a function defined on $\RR$, then the Cauchy integral
$$
 g(z)=\f{1}{2\pi i}\int_{-\infty}^\infty \f{f(t)}{t-z}\,dt
$$
can be written as $(\cP(I+i\Ht)f)(z)$, where $\Ht$ is the Hilbert transform (convolution with $(\pi t)^{-1}$
understood in the principal value sense).
Since the operator $\Ht$ is bounded in $L^p(\RR)$ for $1<p<\infty$, we deduce from (\ref{LpP}):
$$
 \|\cM_g^{(0,\pi)}\|_p\leq K_2(p)\,\left(1+\|\Ht\|_{L^p}\right)\,\|f\|_p,\qquad 1<p<\infty.
$$
For $f$ supported on $[0,\infty)$, the function $g$ is analytic in $\CC^*$ and the estimate identical to
the one just obtained holds true for $\|\cM_g^{(\pi,2\pi)}\|_p$.
Hence (\ref{LpS}) follows with
$$
 K_3(p)=2 \,K_2(p)\,\left(1+\|\Ht\|_{L^p}\right).
$$

\bigskip\noindent
{\bf 5.}\@ To prove Theorem \ref{mLthm},
we refer to the following 
result by the author and A.E.~Merzon \cite[\S~2.1]{MS}, \cite{SM-CMA10}.

Consider the restriction of the Laplace transform $\Lt f$ to a ray $\arg z=\theta$, where
$|\theta|<\pi/2$:
\beq{Ltray}
 \Lta f(\rho)=\int_0^\infty f(t)e^{-t\rho e^{i\theta}}\,dt.
\eeq
Thus, $\Lta$ is an operator that sends a function defined on $\RR_+$ to a function defined on $\RR_+$.
The result mentioned extends the Hausdorff-Young theorem (which formally corresponds to the limiting cases
$\theta=\pm\pi/2$):
\beq{Lpray}
\|\Lta f\|_{p'}\leq K_5\|f\|_p, \qquad 1\leq p\leq 2,
\eeq
where $K_5$ depends on $p$, but is independent of $f$ and $\theta$. We emphasize that the estimate
is uniform with respect to $\theta$.

If $-\pi/2<\theta_1<\theta_2<\pi/2$ and $f(t)$ is a simple (finitely-supported, with finitely many values)
function defined on $\RR_+$, then $\Lt f(z)=O(z^{-1})$ uniformly in $C_{(\theta_1,\theta_2)}$ and
the Cauchy representation readily follows:
$$
 \Lt f(z)=\f{1}{2\pi i}\left(\;\int\limits_{\arg \zeta=\theta_2}-\int\limits_{\arg\zeta=\theta_1}\right)
 \f{\Lt f(\zeta)}{\zeta-z}\,d\zeta,
$$
provided that $\theta_1<\arg z<\theta_2$.

Suppose $1<p\leq 2$ (the case $p=1$ is trivial).
Using the estimate (\ref{Lpray}) on both rays of integration and combining it with Theorem \ref{mSthm},
we conclude that the inequality
$$
\|\cM_{\Lt f}^{(\theta_1,\theta_2)}(\rho)\|_{p'} \leq C\,\|f\|_p
$$
holds for all simple functions $f$, with $C$ independent of $f$ and of $\theta_1, \theta_2$.
Taking the supremum over $(\theta_1,\theta_2)\subset(-\pi/2,\pi/2)$, we obtain (\ref{LpL})
for simple functions $f$; the general case follows by density.
(The operator $\mL$ is nonlinear; however, we can fix some measurable
function $\theta(\rho):\,\RR_+\to (-\pi/2,\pi/2)$, then consider the operator 
$f(t)\mapsto \Lt f(\rho e^{i\theta(\rho)})$;
this operator is linear and the standard density argument can be applied; finally we take supremum
over all functions $\theta(\rho)$.)
 

\bigskip\noindent
{\bf 6.}\@ 
 Finally, we prove Theorem \ref{mwPthm}.
Without loss of generality we may assume $f\geq 0$.

Suppose for simplicity that the supremum in (\ref{mpt}) is attained, specifically ---
at $\theta=\theta_*(R)$.
(The set of all functions for which this is true is dense in $L^1(\RR)$: it contains, say, all simple functions.)
We will estimate the measure of the set
$$
E_{\mP f}(\lambda)\cap\{R\,:\,\theta_*(R)\in(0,\pi/2]\}.
$$
An identical estimate will hold true for
$|E_{\mP f}(\lambda)\cap\{R\,:\,\theta_*(R)\in[\pi/2,\pi)\}|$.

Denote $x_*(R)=R\cos\theta_*(R)$, $y_*(R)=R\sin\theta_*(R)$, and $\delta=\delta(R)=R-x_*(R)$.
By our assumption, $\delta\geq 0$.
Let us split the Poisson kernel as follows:
$$
 P(t,y)=\pi^{-1}\f{y}{y^2+t^2}=P_1(t,y)+P_2(t,y),
$$
where
$$
\ba{l}
\dst
 P_1(t,y)=\min(P(t,y), P(\delta,y)),
\\[3ex]\dst
 P_2(t,y)=\left\{\ba{l} 0, \quad |t|\geq \delta,\\[1ex]
                        P(t,y)-P(\delta,y), \quad |t|<\delta.
\ea\right.
\ea
$$
Denote
$$
 g_k(x,y)=\int P_k(x-t,y) f(t)\,dt, \qquad k=1,2.
$$
(The dependence of $g_k$ on $R$ through the dependence $\delta(R)$ is implicit here.)

\medskip
We will obtain separate estimates of the measure of ``large value sets'' for the functions $g_1$
and $g_2$. Below, $\lambda_1$ and $\lambda_2$ are some arbitrarily chosen positive numbers.
After the separate estimates are obtained, we will specify $\lambda_1$ and $\lambda_2$
in terms of the threshold value $\lambda$ from (\ref{wL1P}).

\medskip
\noindent
\unl{Weak $L^1$ estimate for $g_1(x_*(R), y_*(R))$}

\smallskip
The Poisson kernel is a convex combination of the normalized characteristic functions
of centered segments, and $P_1$ is a sub-convex combination of those corresponding to the
intervals containing the point $\delta$. 
Here is the formal statement.

\begin{lem}
Let the function $\phi_y(a)$, $a>\delta$ be defined as
$$
 \ph_y(a)=-2a\f{d}{da}P_1(a,y).
$$
Then $\ph_y(a)>0$ on $(\delta,\infty)$, $\int_{\delta}^\infty \ph_y(a)\,da<1$,
and
\beq{decomp1}
 P_1(t,y)=\int_{\delta}^\infty \f{1}{2a}\chi_{[-a,a]}(t)\, \ph_y(a)\,da.
\eeq
\end{lem}

\pf
The inequality
$\ph_y(a)>0$ is obvious, since $P_1(t,y)$ is monotonely decreasing in $t$ when $t>\delta$.
Then,
$$
\ba{c}\dst
\int_{\delta}^\infty \ph_y(a)\,da=-2\int_\delta^\infty a\,d_a P_1(a,y)
\\[3ex]\dst
=-2a P_1(a,y)|_{a=\delta}^\infty+2\int_\delta^\infty P_1(a,y)\,da
\\[3ex]\dst
=2\delta P_1(\delta,y)+2\int_\delta^\infty P_1(a,y)\,da
\\[3ex]\dst
=\int_{-\infty}^\infty P_1(a,y)\,da
<\int_{-\infty}^\infty P(a,y)\,da=1.
\ea
$$

Let us verify formula (\ref{decomp1}).
Since the functions $P_1(t,y)$ and $\chi_{[-a,a]}(t)$ are even
(as functions of $t$), we may assume that $t>0$. If $t<\delta$, then $\chi_{[-a,a]}(t)=1$ for all
$a\geq\delta$ and the right-hand side of (\ref{decomp1}) becomes
$$
\int_{\delta}^\infty \f{1}{2a}\, \ph_y(a)\,da=-\int_{\delta}^\infty \f{d}{da}P_1(a,y)\,da=P_1(\delta,y).
$$
If $t\geq\delta$, then the right-hand side of (\ref{decomp1}) becomes
$$
\int_{a>t}\f{1}{2a}\, \ph_y(a)\,da=-\int_{t}^\infty \f{d}{da}P_1(a,y)\,da=P_1(t,y),
$$
as required.
\eop

\medskip
If $g_1(x,y)>\lambda_1$, then, due to Lemma 1, there exists $a>\delta$ such that
$$
\f{1}{2a}\int_{x-a}^{x+a} f(t)\,dt>\lambda_1.
$$
With $x=x_*(R)$, $y=y_*(R)$, $\delta=\delta(R)$ we get
(since $a>\delta$)
$$
\sup_{[u,v]\ni R}\;\f{1}{v-u}\int_{u}^v f(t)\,dt
\;\geq\,
\f{1}{2a}\int_{R-(\delta+a)}^{R+(a-\delta)} f(t)\,dt>\lambda_1.
$$
According to the non-centered Hardy-Littlewood maximal theorem, the set of all $R$ for which such
an inequality holds, has measure not exceeding
\beq{estmu1}
\mu_1= \f{C_1\|f\|_1}{\lambda_1},
\eeq
with some universal constant $C_1$.
%
 Note that the conclusion does not involve $\delta$. 

\medskip
\noindent
\unl{Weak $L^1$ estimate for $g_2(x_*(R), y_*(R))$}

\smallskip
We have
$$
 P_2(t,y)\leq P_2(0,y)=\f{y}{\pi}\left(\f{1}{y^2}-\f{1}{y^2+\delta^2}\right)=
 \f{\delta^2}{\pi y(\delta^2+y^2)}.
$$
Therefore, if
\beq{badR2}
 \int P_2(x_*-t,y_*) f(t)\,dt> \lambda_2,
\eeq
then
$$
 \f{\delta^2}{\pi y_*(\delta^2+y_*^2)}>\f{\lambda_2}{\|f\|_1}.
$$
Now,
\beq{R2xy}
 R^2=x_*^2+y_*^2=(R-\delta)^2+y_*^2=R^2-2\delta R+\delta^2+y_*^2,
\eeq
so
$$
 y_*^2+\delta^2=2R\delta.
$$
Thus, the inequality (\ref{badR2}) implies
$$
 \f{\delta}{2\pi y_*R}>\f{\lambda_2}{\|f\|_1},
$$
hence
$$
 R<\f{\delta\|f\|_1}{2\pi y_*\lambda_2}.
$$
Note that $x_*>0$, so $0<R-\delta$, hence $\delta<2R-\delta$. It follows by (\ref{R2xy}) that
$$
 \left(\f{\delta}{y_*}\right)^2=\f{\delta^2}{\delta(2R-\delta)}<1,
$$
and
$$
 R<\f{\|f\|_1}{2\pi \lambda_2}.
$$
Thus the measure of the set of all such $R$ for which $g_2(x_*(R), y_*(R)>\lambda_2$ does not exceed
$$
 \mu_2=\f{\|f\|_1}{2\pi \lambda_2}.
$$
Again, there is no mentioning of $\delta$ in the final estimate.

To finish the proof of Theorem~\ref{mwPthm}, put $\lambda_1=\lambda_2=\lambda/2$. If $\lambda<g(x_*(R),y_*(R))=g_1(x_*,y_*)+g_2(x_*,y_*)$, then at least one of the
two inequalities holds: $g_1(x_*,y_*)>\lambda_1$ or $g_2(x_*,y_*)>\lambda_2$. Thus
$$
 E_g(\lambda)\subset E_{g_1}(\lambda_1)\cup E_{g_2}(\lambda_2),
$$
so
$$
 \mu_g(\lambda)\leq \mu_1+\mu_2\leq \f{C_1\|f\|_1}{\lambda/2}+\f{1}{2\pi}\,\f{\|f\|_1}{\lambda/2}
 =C\f{\|f\|_1}{\lambda},
$$
with
$C=2C_1+\pi^{-1}$.

\eop

\bigskip\noindent
{\bf 7}. Here are some remarks, open questions, and directions of further work.

\smallskip
1) In Theorem~\ref{mwPthm}, the domain $L^1(\RR)$ of the operator $\mP$ can be extended to the space of measures $\mathcal{M}=C(\RR)^*$.

\smallskip
2) The Poisson convolution operator possesses a natural translation symmetry, which is not respected by the angular maximal function. Translating the
origin of the radial rays in the definition of the angular maximal function leads to trivial generazations of
Theorems~\ref{mwPthm}--\ref{mLthm}. But aren't there more interesting implications --- or a fair translation-invariant version? 

\smallskip
3) In Theorem~\ref{mSthm}, the case $p=1$ is excluded. I expect that, similarly to Theorem~\ref{mwPthm}, the maximal angular Stieltjes transform is of weak type (1,1). To prove this, a method bypassing the Hilbert transform is needed.

\smallskip
4) The technique of proof of Theorem~\ref{mwPthm} here is
similar to that in \cite[second proof of Th.~12]{MS}. 
What property of measure can be axiomatized in order to
put the two theorems in a common framework? A similar question asked about the technique of proof in \cite{SM-CMA10}
has led to the notion of well-projected measures in \cite{MS}.

\smallskip
5) A nice concrete reason for seeking a generalization
mentioned in item (4) is the following anticipated discrete analog of Theorem~\ref{mHilb},
which can also be viewed as a ``maximal interpolation'' between two
classical Hilbert inequalities \cite[8.12 (294)]{HLP}.
\\[0.3ex]
 {\em 
Given two sequences $\{\theta_n\}$ and $\{\theta'_n\}$
with ranges in $(0,\pi)$, 
define the infinite matrix 
$$
A_{mn}=\frac{1-\delta_{mn}}{me^{i\theta_m}-ne^{-i\theta'_n}},
\quad m,n\geq 1.
$$ 
The family of operators with such matrices
is uniformly bounded in $l^p$.}

\smallskip
6) The results presented here are one-dimensional. To what extent are multidimensional generalizations straightforward
and what new interesting subtleties will the multidimensional case bring about?


\begin{thebibliography}{2}

\bibitem{Grafakos}
L.~Grafakos, {\it Classical and Modern Fourier Analysis}, Pearson Education Inc., 2004.

\bibitem{HLP}
G.H.~Hardy, J.E.~Littlewood, G.~Polya, {\em Inequalities}. Cambridge Univ.\ Press, 1934.

\bibitem{MS}
A.~Merzon and S.~Sadov, Hausdorff-Young type theorems for the Laplace transform restricted to
a ray or to a curve in the complex plane. E-print, {\tt http://arxiv.org/math.CA/1109.6085}

\bibitem{SM-CMA10}
S.~Sadov and A.~Merzon,
$L^2$-estimates for the Laplace transform along a family of hyperbolas in the right half-plane,
in: Proceedings of ``Analysis, Mathematical Physics and Applications''
(Ixtapa, Mexico, March 1--5, 2010), 
{\em Comm.\ in Math.\ Analysis},
Conference 03 (2011), 204--208.

\bibitem{Stein} E.~M.~Stein and G.~Weiss, {\em Introduction to Fourier analysis on Euclidean spaces}, Princeton Univ.~Press, Princeton, NJ, 1974.

\end{thebibliography}
\end{document}